%
%
\input amstex
\documentstyle{amsppt}
\loadbold

\def\<{\left<}
\def\>{\right>}

\def\qedd{
\hfill
\vrule height4pt width3pt depth2pt
\vskip .5cm
}
\nologo

\magnification=\magstephalf

\topmatter

\title
The Structure of Spin Systems
\endtitle

\author William Arveson\\
and\\Geoffrey Price
\endauthor

\abstract
A spin system is a sequence of self-adjoint unitary operators 
$U_1,U_2,\dots$ acting on a Hilbert space $H$ which either commute or 
anticommute, $U_iU_j=\pm U_jU_i$ for all $i,j$; it is
is called irreducible when $\{U_1,U_2,\dots\}$ is an irreducible 
set of operators.  
There is a unique infinite 
matrix $(c_{ij})$ with $0,1$ entries satisfying 
$$
U_iU_j=(-1)^{c_{ij}}U_jU_i, \qquad i,j=1,2,\dots.  
$$
Every matrix $(c_{ij})$ with  $0,1$  entries
satisfying $c_{ij}=c_{ji}$ and $c_{ii}=0$ arises from 
a nontrivial 
irreducible spin system, and there are 
uncountably many such matrices.  

Infinite dimensional irreducible representations 
exist when the commutation matrix
$(c_{ij})$ is of ``infinite rank".  In 
such cases we show that the $C^*$-algebra generated 
by an irreducible spin system 
is the CAR algebra, an infinite tensor product of copies 
of $M_2(\Bbb C)$, and we classify the irreducible spin 
systems associated with a given 
matrix $(c_{ij})$ up to 
approximate unitary equivalence.  
That follows from a structural 
result.  The $C^*$-algebra generated by 
the universal spin system $u_1,u_2,\dots$ of $(c_{ij})$ 
decomposes into a tensor product $C(X)\otimes\Cal A$, where 
$X$ is a Cantor set (possibly finite) and 
$\Cal A$ is either the CAR algebra 
or a finite tensor product of copies of $M_2(\Bbb C)$.  

The Bratteli diagram technology of AF algebras 
is not well suited to  
spin systems.  Instead, we work out elementary 
properties of the $\Bbb Z_2$-valued ``symplectic" form 
$$
\omega(x,y) =\sum_{p,q=1}^\infty c_{pq}x_qy_p, 
$$
$x,y$ ranging over the 
free infninite dimensional 
vector space over the Galois field $\Bbb Z_2$,   
and show that one can read off the structure of $C(X)\otimes\Cal A$ 
from properties of $\omega$.  
\endabstract

\affil Department of Mathematics\\
University of California\\Berkeley CA 94720, USA\\\\
Department of Mathematics\\U. S. Naval Academy\\
Annapolis, MD 21402, USA
\endaffil

\rightheadtext{Spin Systems}
\endtopmatter

\thanks
The authors acknowledge support from 
NSF grants DMS-9802474 and DMS-9706441.  
\endthanks

\endtopmatter
%

\document

\subhead{1.  Introduction}
\endsubhead

A spin system is a sequence $u_1,u_2,\dots$ of self-adjoint 
unitary elements of some unital $C^*$-algebra which 
commute up to phase in the sense that 
$$
u_iu_j=\lambda_{ij}u_ku_j, \qquad i,j =1,2,\dots
$$
where the $\lambda_{ij}$ are complex numbers.  Since 
$u_iu_ju_i^{-1}=\lambda_{ij}u_j$ and $u_j^2=\bold 1$, it 
follows that each $\lambda_{ij}$ is $-1$ or $+1$.  Thus there 
is a unique matrix of zeros and ones $c_{ij}$ such that the 
commutation relations become 
$$
u_iu_j=(-1)^{c_{ij}}u_ju_i, \qquad i,j=1,2,\dots.  \tag{1.1}
$$
The matrix $(c_{ij})$ is symmetric and 
has zeros along the main diagonal.  
A concrete spin system $U_1,U_2,\dots\subseteq\Cal B(H)$ is 
said to be irreducible when $\{U_1,U_2,\dots\}$ is an 
irreducible set of operators.  The purpose of this paper 
is to determine the structure of the $C^*$-algebra generated 
by an irreducible spin system associated with a given 
$0$-$1$ matrix $(c_{ij})$, and to classify such spin 
systems up to ``approximate" unitary equivalence (Theorem 
C, section 4).   

\remark{Quantum Spin Systems}
Spin systems arise naturally in several contexts, including 
the theory of quantum spin systems (\cite{BR}, section 6.2), 
and in the theory of quantum computing (especially,
systems involving a large or infinite number of qubits).  
For example, suppose we are given a mutually 
commuting sequence $\theta_1, \theta_2, \dots$ of involutive 
$*$-automorphisms of $\Cal B(H)$, i.e., $\theta_j^2={\text{id}}$, 
$\theta_k\theta_j=\theta_j\theta_k$ for all $j,k=1,2,\dots$
(one can imagine that $\theta_k$ represents reversing 
the state of a two-valued quantum observable
located at the $k$th site).  
For each $k$ one can find a unitary operator $U_k$ such that 
$\theta_k(A)=U_kAU_k^{-1}$, $A\in\Cal B(H)$, 
and by replacing $U_k$ with 
$\lambda U_k$ for an appropriate scalar $\lambda$ if 
necessary, we can arrange that $U_k^2=\bold 1$.  Since 
$\theta_i\theta_j=\theta_j\theta_i$ it follows that $U_i$ 
and $U_j$ must commute up to phase, hence there is a 
unique number $c_{ij}\in\{0,1\}$ such that (1.1) is 
satisfied.  The matrix $C=(c_{ij})$ does not depend on 
the choices made and is therefore an invariant 
attached to the original sequence of automorphisms 
$\bar \theta=(\theta_1,\theta_2,\dots)$.  
The sequence $\bar\theta$ is {\it ergodic}
in the sense that  its fixed algebra is $\Bbb C\cdot\bold 1$ if, and 
only if, the spin systems $\bar U=(U_1, U_2,\dots)$ 
associated with it are irreducible.  
\endremark

\remark{Remarks on rank}
Consider the commutation matrix $(c_{ij})$ associated with 
a spin system (1.1).  If all coefficients $c_{ij}$ vanish then 
$C^*(u_1,u_2,\dots)$ is commutative.  More generally, 
$C^*(u_1,u_2,\dots)$ degenerates whenever $(c_{ij})$ is of 
finite rank, where the rank is defined as follows.  
Considering $\Bbb Z_2=\{0,1\}$ as the two-element 
Galois field we may consider vector spaces over $\Bbb Z_2$, 
and in particular we can form the free infinite
dimensional  vector space 
$\Gamma=\Bbb Z_2\oplus\Bbb Z_2\oplus\dots$ over $\Bbb Z_2$.  
Elements of $\Gamma$ are sequences $x=(x_1,x_2,\dots)$,
$x_k\in\Bbb Z_2$, which vanish eventually.  The 
dual of $\Gamma$ is identified 
with the vector space $\Bbb Z_2^\infty$ 
of all sequences $y=(y_1, y_2,\dots)$, $y_k\in\Bbb Z_2$.  
The commutation matrix $(c_{ij})$ gives rise to a linear 
operator $C: \Gamma\to\Bbb Z_2^\infty$ by way of 
$(Cx)_k=\sum_{j=1}^\infty c_{kj}x_j$, $k=1,2,\dots$.  
The {\it rank} of the matrix is defined by 
$$
{\text{rank}}\,(c_{ij}) = \dim(C\Gamma).  
$$
Actually, what we have defined is the {\it column rank} 
of the matrix $(c_{ij})$, but because $(c_{ij})$ is a 
symmetric matrix its column and row ranks are the same.  
We will see below that the rank is finite iff 
the center of $C^*(u_1,u_2,\dots)$ is of finite codimension
in $C^*(u_1,u_2,\dots)$ iff every irreducible spin system 
satisying (1.1) acts on a finite dimensional Hilbert 
space.  
Thus we are primarily concerned with the nondegenerate 
cases in which the commutation matrix $(c_{ij})$ is 
of infinite rank.  
\endremark

\remark{Remarks on Existence and Universality}
A sequence $u_1,u_2,\dots$ of unitary operators satisfying a 
given set of 
noncommutative equations $f_k(u_1, u_2,\dots)=0$, $k=1,2,\dots$
(we leave the precise nature of the noncommutative polynomials 
$f_k$ unspecified) is said to be universal if
every  sequence $U_1,U_2,\dots\in\Cal B(H)$ of 
concrete unitary operators that 
satisfies the equations can be obtained from it via 
a representation $\pi:C^*(u_1,u_2,\dots)\to\Cal B(H)$ 
such that $\pi(u_k)=U_k$, $k=1,2,\dots$.  

Of course, for  bad choices of $f_k$ (such as 
$f_k(x_1,x_2,\dots)=x_k$)  there may be no unitary 
solutions to the set of equations except 
on the trivial Hilbert space $H=\{0\}$.  But in all 
cases there is a universal solution...the direct 
sum of all concrete unitary
solutions.   Any two universal solutions $(u_1,u_2,\dots)$ and 
$(v_1,v_2,\dots)$ are equivalent in the sense that there 
is a unique $*$-isomorphism $\theta: C^*(u_1,u_2,\dots)\to
C^*(v_1,v_2,\dots)$ satisfying $\theta(u_k)=v_k$ for 
every $k$. Thus the $C^*$-algebra generated by a universal 
sequence of solutions to the given set $S$ of equations is 
uniquely determined by $S$.  
\endremark 

Given an arbitrary matrix $C=(c_{ij})$ 
of zeros and ones satisfying the consistency requirements 
$c_{ij}=c_{ji}$ and $c_{jj}=0$ for all $i,j=1,2,\dots$, we 
consider the $C^*$-algebra $\Cal A_C=C^*(u_1,u_2,\dots)$  
generated by a {\it universal} spin system satisfying (1.1).  
The set of distinct matrices $(c_{ij})$ satisfying 
these conditions is of cardinality $2^{\aleph_0}$, 
and each of them 
is associated with a nontrivial spin system (1.1)  
(see Proposition (1.1)).    We determine the structure 
of these $C^*$-algebras $\Cal A_C$ in Theorem B, section 3.

\remark{Spin Systems in Characteristic $p$}
We have found it helpful, even simplifying,
to consider the natural generalization of spin systems to 
characteristic $p$ where $p$ is an arbitrary prime.  
By a spin system in characteristic $p$ we mean a sequence of unitary 
operators $u_1, u_2,\dots$ which are $p$th roots of unity in the 
sense that $u_j^p=\bold 1$ for every $j$, and which satisfy 
commutation relations of the form 
$$
u_iu_j=\zeta^{c_{ij}}u_ju_i, \qquad i,j=1,2,\dots,\tag{1.2}
$$
where $\zeta=e^{2\pi i/p}$, and where 
$c_{ij}\in\{0,1,\dots,p-1\}=\Bbb Z_p$.  The matrix $(c_{ij})$ 
is uniquely determined by (1.2).  If we regard $\Bbb Z_p$ as 
a finite field in the usual way, then the matrix 
is {\it skew}-symmetric in that $c_{ij}=-c_{ji}$ for every 
$i,j = 1,2,\dots$.  

The reason for considering the cases 
$p>2$ can be clearly seen 
when one specializes the previous paragraph to
$p=2$.   Indeed, 
a skew-symmetric matrix 
over the two-element field $\Bbb Z_2$ is the same as a 
symmetric matrix with zeros along the main diagonal.  We 
found that viewing $(c_{ij})$ as a
skew-symmetric matrix  led in the right direction, 
whereas viewing it as a symmetric matrix with 
zeros along the diagonal led nowhere.  Thus the 
case $p=2$ can be misleading, and for that reason 
we consider the 
more general case of spin systems (1.2) in 
characteristic $p$.  
\endremark

Fixing a prime $p$, suppose we are given a skew-symmetric 
matrix $(c_{ij})$ of elements of 
the Galois field $\{0,1,\dots,p-1\}=\Bbb Z_p$.  
Since $\Bbb Z_p$ is a field, we can form the 
free infinite dimensional vector space 
$\Gamma$ over $\Bbb Z_p$; elements of $\Gamma$ are sequences 
$x=(x_1,x_2,\dots)$ of elements of $\Bbb Z_p$ satisfying 
$x_k=0$ for all but a finite number of $k$.  The coefficients 
$c_{ij}$ give rise to a bilinear form 
$\omega: \Gamma\times\Gamma\to\Bbb Z_p$ by way of 
$$
\omega(x,y) = \sum_{i,j=1}^\infty c_{ij}x_jy_i, 
\qquad x,y\in\Gamma.  \tag{1.3}
$$
This bilinear form is {\bf skew-symmetric} in that it 
satisfies $\omega(x,y)=-\omega(y,x)$ for 
all $x,y\in \Gamma$, and it will occupy 
a central position throughout the sequel.  The structure 
of such forms is described in Theorem A, section 2, 
and its Corollary.    

Consider now the 
$C^*$-algebra $\Cal A$ generated by 
a sequence of unitary elements $u_1, u_2,\dots$ satisfying 
$u_k^p=\bold 1$ and the commutation relations (1.2).  
A {\bf word} is a finite product of elements from $\{u_1,u_2,\dots\}$, 
and it is convenient to regard the identity $\bold 1$ as the 
empty word.  The set of linear combinations of words is 
a dense $*$-subalgebra of $\Cal A$ which contains $\bold 1$.  
Using the commutation relations (1.2), every word can be 
written in the form 
$\lambda u_1^{n_1}u_2^{n_2}\dots u_r^{n_r}$ where 
$\lambda$ is a complex scalar.  Thus we may use the 
elements of $\Gamma$ to parameterize a spanning set of 
words as follows,
$$
w_x=u_1^{x_1}u_2^{x_2}\dots, \qquad x=(x_1,x_2,\dots)\in\Gamma, 
$$
and one finds that 
$$
w_xw_y=\zeta^{\omega(x,y)}w_yw_x, \qquad x,y\in \Gamma,  \tag{1.4}
$$
where $\zeta=e^{2\pi i/p}$ and 
$\omega: \Gamma\times\Gamma\to\Bbb Z_p$ is the 
bilinear form (1.3).   

We will occasionally make use of a second bilinear form 
$Q: \Gamma\times\Gamma\to\Bbb Z_p$,
$$
Q(x,y) = \sum_{1\leq i<j}c_{ij}x_jy_i, \qquad x,y\in \Gamma.  \tag{1.5}
$$
$Q$ is related to $\omega$ by 
$\omega(x,y)=Q(x,y)-Q(y,x)$, and it 
obeys the ``Weyl" relations
$$
w_xw_y=\zeta^{Q(x,y)}w_{x+y}, \qquad x,y\in\Gamma.  \tag{1.6}
$$ 

We conclude the introduction with an observation 
about the existence of 
solutions of (1.2) for arbitrary coefficient matrices $(c_{ij})$.  
For $p=2$, this generalizes the examples of finite 
dimensional spin systems described in \cite{Bi}.  

\proclaim{Proposition 1.1}
Let $p=2,3,\dots$ be a prime and let $(c_{ij})$ be an 
arbitrary  skew-symmetric 
matrix over the Galois field 
$\Bbb Z_p=\{0,1,\dots,p-1\}$.  There 
a Hilbert space $H\neq\{0\}$ 
and a sequence of unitary operators 
$U_1,U_2,\dots\in\Cal B(H)$ such that $U_k^p=\bold 1$ 
and $U_jU_k=\zeta^{c_{jk}}U_kU_j$ for every $j,k=1,2,\dots$, 
where $\zeta=e^{2\pi i/p}$.  
\endproclaim

\demo{proof}
Regarding $\Bbb Z_p$ as an additive abelian group, consider the 
unitary operators $S$, $V$ defined on 
the $p$-dimensional Hilbert space
$\ell^2(\Bbb Z_p)$ by 
$$
Sf(k)=f(k+1), \quad Vf(k)=\zeta^kf(k), \qquad f\in \ell^2(\Bbb Z_p), 
\quad k\in\Bbb Z_p.  
$$
We have $S^p=V^p=\bold 1$, 
$SV=\zeta VS$, and in fact 
$SV^k=\zeta^kV^kS$ for all $k\in \Bbb Z$.  

Consider the $L^2$-space of the compact abelian group 
$G=\Bbb Z_p\times\Bbb Z_p\times\dots$.  We can realize 
$L^2(G)$ as the infinite tensor product of copies of 
$\ell^2(\Bbb Z_p)$ along the stabilizing vector $u\in\ell^2(\Bbb Z_p)$
where $u$ is the constant function $u(k)=1$, $k\in\Bbb Z_p$.  Thus 
for any finite sequence $A_1,\dots,A_r$ of operators on 
$\ell^2(\Bbb Z_p)$ we can form the operator 
$$
A_1\otimes\dots\otimes A_r\otimes\bold 1\otimes\bold 1\otimes\dots
\in\Cal B(L^2(G)).  
$$
The unitary operators $U_1,U_2,\dots$ are defined on $L^2(G)$ 
in terms of the given 
coefficients $c_{ij}$ as follows; 
$U_1=S\otimes\bold 1\otimes\bold 1\otimes\dots$ and for 
$k= 2,3,\dots$
$$
U_k=V^{c_{1k}}\otimes\dots\otimes V^{c_{k-1 k}}\otimes 
S\otimes\bold 1\otimes\bold
1\otimes\dots.  
$$
One can verify that $U_k^p=\bold 1$, and 
$U_jU_k=\zeta^{c_{jk}}U_kU_j$ for $1\leq k<j$.\qedd
\enddemo

As the preceding remarks on Quantum Spin Systems show, 
commutation relations of the form $uv=\lambda vu$ where 
$\lambda\in\Bbb T$ arise whenever one 
considers commuting $*$-automorphisms of $\Cal B(H)$, 
and in fact many natural contexts lead 
to projective representations of groups involving
similar commutation relations.  
For example, they are associated with 
ergodic actions of compact groups 
on $C^*$-algebras (see \cite{BE}, \cite{J} 
and references therein).  
Since commutation relations of this type 
are so ubiquitous, 
we have made no effort to compile references 
to the related literature, even for the case of 
spin systems.  
Finally, we point out that the results of this 
paper generalize certain 
results in \cite{PP}, \cite{P1}, \cite{P2} 
which concern spin systems for which the commutation matrix 
depends only on the separation $c_{ij}=f(i-j)$.

\subhead{2.  Symplectic forms in characteristic $p$}
\endsubhead

In this section we work out the symplectic linear algebra that 
underlies the results described above.  While the results 
(and methods) are quite elementary, we have been 
unable to find what we require in the literature;
indeed, fields of characteristic $2$ are 
excluded from most treatments of 
linear algebra.  Thus 
we provide a complete discussion.  
Throughout, $F$ denotes a field, the 
primary cases being the Galois field $F=\Bbb Z_p$ 
of characteristic $p$ where $p$ is any prime including 
$2$.  $\Gamma$ denotes the free infinite dimensional 
vector space over $F$, consisting of all sequences 
$x=(x_1,x_2,\dots)$ of elements $x_k\in F$ satisfying 
$x_k=0$ for all but a finite number of $k$.  The addition 
and scalar multiplication are defined pointwise, 
$$
x+y=(x_1+y_1,x_2+y_2,\dots), \qquad 
\lambda\cdot x=(\lambda x_1,\lambda x_2,\dots)
$$
$\lambda$ being an element of $F$.  A vector space 
$V$ over $F$ is said to be {\it countably generated}
if it contains a sequence $v_1,v_2,\dots$ such that 
every element of $v$ is a finite linear combination 
of elements of $\{v_1,v_2,\dots\}$.  For every countably 
generated vector space $V$ over $F$ there is a linear 
map $L:\Gamma\to V$ such that $V=L\Gamma$.  

We are concerned with skew-symmetric bilinear forms 
$B: \Gamma\times\Gamma\to F$.  
The {\bf kernel} of such a bilinear form is the 
subspace $K=\{x\in\Gamma: B(x,\Gamma)=\{0\}\}$.  
$B$ is called a {\bf symplectic} form when it is 
skew-symmetric and has kernel $\{0\}$, and a 
{\bf symplectic vector space} is a pair 
$(V,B)$ consisting of a countably generated 
vector space $V$ over $F$ and a symplectic 
bilinear form $B: V\times V\to F$.  Two 
symplectic vector spaces $(V,B)$ and 
$(V^\prime,B\prime)$ are {\bf congruent} 
if there is a linear isomorphism $L: V\to V^\prime$ 
satisfying $B^\prime(Lx,Ly)=B(x,y)$ for 
all $x,y\in V$.

Any skew-symmeteric bilinear form defined on a vector 
space $B:V\times V\to F$ gives rise to a symplectic 
vector space as follows.  Letting $K$ be the kernel of $B$, 
$B$ promotes natrually 
to a bilinear form $\omega: V/K\times V/K\to F$,
$$
\omega(x+K,y+K)=B(x,y), \qquad x,y\in V.  \tag{2.1}
$$ 
$(V/K,\omega)$ is a symplectic vector space, 
and it is the trivial symplectic vector space 
only when $B=0$.  

\proclaim{Definition}
The rank of $B$ is the dimension 
of the vector space $V/\ker B$.  
\endproclaim

The rank of $B$ is a nonnegative integer 
or $\infty$.  We will see presently that when 
it is finite it must be an even integer $n=2r$, $r=1,2,\dots$. 

\remark{Further remarks on rank}
Let $C=(c_{ij})$ be the commutation matrix associated 
with the relations (1.1).  We have given a different 
definition of rank 
in the introduction, and we want to point out that 
the rank defined in the introduction 
is the same as the rank of the skew-symmetric 
form $\omega$ associated to it by (1.3).  To see that 
consider the linear map $L: \Gamma\to \Bbb Z_2^\infty$ 
defined by 
$$
Lx = (\omega(u_1,x),\omega(u_2,x),\dots), 
$$
where $u_1, u_2, \dots$ is the usual sequence of 
basis vectors in $\Gamma$, $u_k(j)=\delta_{kj}$.  
Noting that the $k$th component of 
$Lx$ is $\omega(u_k,x)=\sum_j c_{kj}x_j$,
one sees that the range of $L$ is the linear 
span of the columns of $(c_{ij})$ and hence 
its dimension is the rank of the matrix $(c_{ij})$.  
On the other hand, the kernel of $L$ is exactly 
$\ker\omega$, 
so that ${\text{rank}}\,C=\dim L\Gamma = \dim(\Gamma/\ker\omega)$, 
as asserted.  
\endremark

Let $(V,\omega)$ be a symplectic vector space. 
 By a {\bf symplectic basis} for $V$ we mean a pair 
of sequences $e_1, e_2, \dots, f_1, f_2,\dots\in V$ with the
properties
$$
\omega(e_i,e_j)=0, \quad
\omega(f_i,f_j)=0, \quad
\omega(e_i,f_j)=\delta_{ij}, \tag{2.2}
$$
for all $i,j=1,2,\dots$ and which span $V$ in the 
sense that every element of $V$ is a finite linear 
combination of the elements $\{e_i, f_j\}$.  The sequences 
are allowed to be either finite or infinite, but if one 
of them 
is finite then the other is also finite of the 
same length.  A simple argument shows that any finite set 
of $2r$ vectors $e_1,\dots,e_r,f_1,\dots,f_r$ which 
satisfy the relations (2.2) must be linearly independent.  
Thus a symplectic basis for $V$ is a countable basis, 
and in particular $V$ must be countably generated.   

\example{The Standard Examples}
We describe the standard models of symplectic vector 
spaces of dimension $n=2,4,6,\dots,\infty$ over an 
arbitrary field $F$.    
Consider first the case $n=\infty$.   Let 
$\Gamma=F\oplus F\oplus\dots$ be the vector 
space of all infinite 
sequences $x=(x_1,x_2,\dots)$, where 
$x_k\in F$ and $x_k=0$ for all but 
a finite number of $k$.  
The symplectic space $(V_\infty,\omega_\infty)$ is 
defined by $V_\infty=\Gamma\oplus\Gamma$ and 
$$
\omega_\infty((x,y), (x^\prime,y^\prime))=
\sum_{k=1}^\infty y_kx_k^\prime-x_ky_k^\prime.
$$
$(V_\infty,\omega_\infty)$
is a countably generated infinite dimensional symplectic 
vector space, and it has a natural symplectic basis 
$\{e_j,f_k\}$, defined by 
$$
e_k=(u_k,0), \quad f_k=(0,u_k), \qquad k=1,2,\dots
$$
where $u_k$ is the standard unit vector $u_k(j)=\delta_{kj}$.

For $n=2r$ finite, we take $V_n$ to be 
the $2r$ dimensional subspace 
$F^r\oplus F^r\subseteq V_\infty$ 
and define $\omega_n$ by restricting $\omega_\infty$ 
to $V_n$.  
\endexample

The following result implies that any two  
countably generated symplectic vector spaces
of the  same dimension are congruent.  

\proclaim{Theorem A}Let $F$ be a field of arbitrary 
characteristic.
\roster
\item"{A1}"
Every countably generated symplectic vector space 
$(V,\omega)$ over $F$ 
has a symplectic basis.  When the dimension of
$V$ is finite it must be  even, $\dim V=2r$, $r=1,2,\dots$.  
\item"{A2}"
Let $\omega$ be a skew-symmetric bilinear form on a countably 
generated vector space $V$, let $K$ be the kernel of $\omega$ 
and let $L$ be any vector space complement $V=K\oplus L$.  
Then the restriction $\omega_L$ of $\omega$ to $L$ is a symplectic 
form.  If $L^\prime$ is any other complement 
$V=K\oplus L^\prime$, then the symplectic spaces 
$(L,\omega_L)$ and $(L^\prime,\omega_{L^\prime})$ are congruent.  
\endroster
\endproclaim

\demo{proof of (A1)}
We first treat the simple case in which 
$V$ is finite dimensional and nonzero.  
Choose any vector $e_1\neq 0$ in $V$.  By nondegeneracy, 
there is a vector $f_1\in V$ with $\omega(e_1,f_1)=1$.  
In order to continue inductively, we require

\proclaim{Lemma 2.1}
Let $(V,\omega)$ be a finite dimensional symplectic 
vector space, let $S\subseteq V$ be a subspace such 
that the restriction of $\omega$ to $S\times S$ is 
nondegenerate, and 
let $K$ be its symplectic complement
$$
K=\{x\in V: \omega(x,S)=\{0\}\}.  
$$
Then $V=S\oplus K$.  
\endproclaim
\demo{proof}
Obviously, $S\cap K=\{0\}$ because the restriction 
of $\omega$ to $S\times S$ is nondegenerate.  We have 
to show that $V=S+K$, and since the intersection of 
these two spaces is trivial 
it suffices to show that $\dim S+\dim K = \dim V$.  

Assuming $S\neq\{0\}$, let $v_1,\dots,v_r$ be a basis 
for $S$, and consider the linear map 
$L: V\to F^r$ defined by 
$$
Lx=(\omega(x,v_1),\dots,\omega(x,v_r)), \qquad x\in V.  
$$
The kernel of $L$ is $K$, and 
we claim that $LV=F^{r}$.  To prove that
we show that the only linear functional $f: F^{r}\to F$ 
that vanishes on $LV$ is $f=0$.  Indeed, writing 
$$
f(t_1,\dots,t_r)=
\sum_{k=1}^r \lambda_kt_k,
$$
for certain $\lambda_j\in F$, the 
vector $v=\sum_k\lambda_kv_k\in S$ satisfies 
$\omega(x,v)=f(Lx)=0$ for all $x\in V$.  Since $\omega$ 
is nondegenerate we must have $v=0$, hence 
$\lambda_1=\dots=\lambda_r=0$, hence $f=0$.  
We conclude that 
$$
\dim V=\dim{\text{ran}}\,L+\dim\ker L =
\dim F^r+\dim K=\dim S+\dim K,
$$
since $\dim S=r=\dim F^r$.  \qedd
\enddemo

Inductively, suppose we have vectors 
$e_1,\dots,e_r,f_1,\dots,f_r\in V$ which satisfy the 
symplectic requirements (2.2) insofar as they make sense, 
and let $S$ be the subspace of $V$ spanned by 
$\{e_k,f_j:1\leq j,k\leq r\}$.   Since 
$\{e_k,f_j\}$ is a symplectic basis for the restriction 
of $\omega$ to $S\times S$, the latter must be nondegenerate.  
By Lemma 2.1 we have $V=S+K$ where 
$K=\{x\in V: \omega(x,S)=\{0\}\}$.  
Thus we can choose a nonzero vector $e_{r+1}$ in $K$.  
Since $\omega(e_{r+1},S)=\{0\}$ and $V=S+K$, there must 
be a vector $f_{r+1}\in K$ for which $\omega(e_{r+1},f_{r+1})=1$. 
An inductive argument completes the proof in the case 
where $V$ is finite dimensional.  

\remark{Remark}
Notice that the preceding argument implies that in a finite 
dimensional symplectic vector space $(V,\omega)$, any 
set of vectors $e_1,\dots,e_r,f_1,\dots,f_r\in V$, 
which satisfy the relations (2.2), can be enlarged to 
a symplectic basis for $V$.  It also shows that a finite 
dimensional symplectic vector space over an arbitrary 
field has even dimension $2\cdot r$, $r=1,2,\dots$.  
\endremark

\vskip0.1truein

Turning now to the infinite dimensional case, we 
claim that there is an increasing sequence 
of finite dimensional subspaces 
$E_1\subseteq E_2\subseteq\dots\subseteq V$ with 
$\cup_n E_n=V$, such that the restriction of 
$\omega$ to $E_n\times E_n$ is nondegenerate for 
every $n$.  Suppose for the moment that this has 
been established.  The preceding paragraphs show 
that we can find a symplectic basis for $E_1$.  
Since the restriction of $\omega$ to $E_2\times E_2$ 
is a symplectic form on $E_2$, the preceding remark
implies that this symplectic
set can be enlarged  to a symplectic 
basis for $E_2$.  Continuing 
inductively, we obtain an increasing sequence 
of symplectic sets, 
each one being a basis for its corresponding 
linear span $E_n$, $n=1,2,\dots$, and 
their  union is a symplectic basis for $\cup_n E_n=V$.  

Thus we have reduced the proof of (A1) to showing that 
there is such a sequence $E_1\subseteq
E_2\subseteq\dots$.  

\proclaim{Lemma 2.2}
Let $(V,\omega)$ be a symplectic vector space 
and let $E$ be a finite dimensional subspace of $V$.  
Then there is a subspace $E^\prime\supseteq E$ of dimension
at most $2\cdot\dim E$ such that the restriction of 
$\omega$ to $E^\prime\times E^\prime$ is nondegenerate.  
\endproclaim

\demo{proof}
Let $K=\{x\in E: \omega(x,E)=0\}$ be the kernel of the 
restriction of $\omega$ to $E\times E$, 
and let $k_1,\dots,k_r$ be a basis for $K$.  We claim
that there are vectors $\ell_1,\dots,\ell_r\in V$ such that 
$$
\omega(k_i,\ell_j)=\delta_{ij}, \qquad 1\leq i,j\leq r.  \tag{2.3}
$$
To see that, consider the $r$-dimensional vector space 
$F^r=\{(t_1,\dots,t_r): t_i\in F\}$, and consider the linear 
map $L: V\to F^r$ defined by 
$$
L(x)=(\omega(k_1,x),\omega(k_2,x),\dots,\omega(k_r,x)), \qquad x\in V.  
$$
We have to show that $L$ is onto: $L(V)=F^r$.  To prove 
that, we show that the only linear functional $f: F^r\to F$ 
which vanishes on the range of $L$ is the zero 
functional.  Choosing such an $f$, we can write 
$$
f(t_1,\dots,t_r)=\lambda_1t_1+\dots+\lambda_rt_r
$$
for a unique $r$-tuple of scalars $\lambda_k\in F$.  Since 
$f(L(x))=0$ for all $x\in V$ we have 
$$
\omega(\sum_{j=1}^r\lambda_jk_j,x)=
\sum_{j=1}^r\lambda_j\omega(k_j,x)=f(L(x))=0.  
$$
By nondegeneracy, we must have $\sum_j\lambda_jk_j=0$, 
hence $\lambda_1=\dots=\lambda_r=0$ because $k_1,\dots,k_r$
are linearly independent, thus (2.3) is proved.  

Setting $L={\text{span}}\{\ell_1,\dots,\ell_r\}$, 
notice that (2.3) implies that the restriction of 
$\omega$ to $K\times L$ is nondegenerate in the 
sense that for every $k\in K$, 
$$
\omega(k,\ell)=0, \quad {\text{for all }}\ell\in
L\implies k=0,  \tag{2.4}
$$
while for every $\ell\in L$, 
$$
\omega(k,\ell)=0, \quad {\text{for all }}k\in K
\implies \ell=0.  \tag{2.5}
$$

Choose such a set of vectors $\ell_1,\dots,\ell_r\in V$, 
let $L={\text{span}}\{\ell_1,\dots,\ell_r\}$, and define 
$E^\prime = E + L$.  
We show that the restriction of $\omega$ to $E^\prime\times E^\prime$ 
is nondegenerate.  For that, suppose that $z\in E^\prime$ 
has the property that $\omega(z,z^\prime)=0$ for every 
$z^\prime\in E^\prime$.  We can write $z=x+\ell$ where 
$x\in E$ and $\ell\in L$ .   Then 
$$
\omega(z,z^\prime)=\omega(x,z^\prime)+\omega(\ell,z^\prime)=0
$$
for all $z^\prime\in E^\prime$.  Picking $z^\prime\in K$ and noting 
that $\omega(x,K)=\{0\}$ (by definition of $K$), we conclude 
that $\omega(\ell,z^\prime)=0$ for all $z^\prime\in K$.  
Because of (2.5), we conclude that $\ell=0$.  
Hence $\omega(x,E^\prime)=\{0\}$.  Since $x\in E\subseteq E^\prime$
this implies that $x$ is an element of $K$ for which 
$\omega(x,E^\prime)=0$.  By (2.4) this implies $x=0$. \qedd
\enddemo

The proof of (A1) is completed as follows.   
Since $V$ is countably generated 
there is a spanning sequence of nonzero 
vectors $v_1,v_2,\dots\in V$; we will construct an increasing 
sequence $E_n$ of finite dimensional subspaces 
such that $E_n$ contains
$v_1,\dots,v_n$ and  the restriction of $\omega$ to $E_n$ is
nondegenerate.   Since $v_1\neq 0$ and $\omega$ is nondegenerate,
choose any $w\in v$ such that $\omega(v_1,w)=1$, and 
set $E_1={\text{span}}\{v_1,w\}$.  The restriction of 
$\omega$ to $E_1$ is nondegenerate because $\{v_1,w\}$ 
is a symplectic basis.  

Suppose now that we have finite dimensional subspaces 
$E_1\subseteq\dots\subseteq E_n$ such that $E_k$ 
contains $v_1,\dots,v_k$ and the restriction of 
$\omega$ to each $E_k$ is nondegenerate.  Applying 
Lemma 2.2 to the space spanned by $E_n$ and $v_{n+1}$, 
we find a finite dimensional space $E_{n+1}$ containing 
$v_{n+1}$ and $E_n$ such that the restriction of 
$\omega$ to $E_{n+1}\times E_{n+1}$ is nondegenerate.  
An induction completes the proof of (1).

In order to prove (A2), consider the natural 
symplectic space $(V/K,\omega)$ described above.  
We claim that for every subspace $L$ 
of $V$ satisfying $L\cap K=\{0\}$ and $L+K=V$,
the symplectic spaces 
$(L,\omega_L)$ and $(V/K,\omega)$ are congruent; 
i.e., there is a linear  isomorphism 
$T: L\to V/K$ such that 
$$
\omega(Tx,Ty)=B(x,y)=\omega_L(x,y), \qquad x,y\in L,  \tag{2.6}
$$
where $\omega_L$ is the
restriction of $B$ to $L\times L$.   To see that, 
define $Tx=x+K$, $x\in L$.  $T$ is a linear isomorphism 
because $L$ is a complement of $K$, and (2.6) follows because 
for any $x,y\in V$ we have 
$\omega(x+K,y+K)=B(x,y)$ by definition of $\omega$, so when 
$x,y\in L$ we have (2.6).  

For any other subspace $L^\prime$ with 
$V=K\oplus L^\prime$, $(L^\prime,\omega_{L^\prime})$ is 
also congruent to $(V/K,\omega)$, hence it is congruent 
to $(L,\omega_L)$.  
\qedd
\enddemo

\proclaim{Corollary}
Any two countably generated symplectic vector spaces 
of the same dimension $n=2r$, $r=1,2,\dots,\infty$ are 
congruent. 
\endproclaim

\demo{proof}
Let $(V,\omega)$ be a symplectic vector space of 
dimension $n=2r$, $r=1,2,\dots,\infty$.  
By Theorem A, we can find a (finite or infinite) symplectic 
basis $\{e_k,f_j\}$ for $V$, and once we have that there 
is an obvious way to transform $(V,\omega)$ congruently to 
the standard example $(V_n,\omega_n)$.  \qedd
\enddemo

\example{Examples of commutation matrices}
The above results have concrete implications about 
how to exhibit sequences of unitary operators 
that generate the infinite dimensional 
CAR algebra; they also provide a systematic  
method for generating all possible 
skew-symmetric matrices $C=(c_{ij})$ with entries in 
$\Bbb Z_2$ which are {\it nondegenerate} in the sense 
that their associated bilinear forms 
$$
\omega_C(x,y)=\sum_{i,j=1}^\infty c_{ij}x_jy_i, 
\qquad x,y\in\Gamma  \tag{2.7}
$$
have trivial kernel.  We abuse 
terminology slightly by 
calling such a matrix $C$ {\it symplectic}.  
Starting with any countably infinite symplectic vector space 
$(V,\omega)$ over $\Bbb Z_2$, such as the standard 
example $(V_\infty,\omega_\infty)$ described above,  
let $v_1,v_2,\dots$ be any basis for $V$ and define 
$$
c_{ij}=\omega(v_i,v_j), \qquad i,j=1,2,\dots.  
$$
One verifies directly that $C=(c_{ij})$ is a 
symplectic matrix.  Moreover, the Corollary 
of Theorem A implies that every symplectic matrix 
arises in this way from some basis 
$v_1, v_2,\dots$ for $V$.  

One can view this construction in more 
concrete operator-theoretic 
terms by making use of the standard 
self-adjoint generators of the 
CAR algebra as follows.  Consider the 
Clifford algebra $\Cal C$ generated by an infinite sequence 
$W_1, W_2,\dots$ of unitary operators satisfying 
$$
W_iW_j + W_jW_i = 2\delta_{ij}\bold 1, 
\qquad i,j=1,2,\dots.  
$$
Since $W_i$ and $W_j$ anticommute when 
$i\neq j$, the commutation matrix $A=(a_{ij})$
associated with a Clifford sequence is  
$$
a_{ij} = 
\cases
1, \quad  &i\neq j,\\
0, \quad  &i=j,  
\endcases
$$
and its associated form is 
$$
\omega_A(x,y) =\sum_{p\neq q}x_qy_p=
(\sum_kx_k)(\sum_ky_k)-\sum_kx_ky_k.  
$$
One verifies easily that $\omega_A$ is nondegenerate.   
Choosing an arbitrary basis $v_1, v_2,\dots$ 
for $\Gamma$, we obtain the most general symplectic matrix 
$C=(c_{ij})$ as follows
$$
c_{ij}=\omega_A(v_i,v_j)=\sum_{p\neq q} v_i(q)v_j(p).  
\tag{2.8}
$$
Each element $v_k$ in this basis is associated with 
a word in the original sequence $(W_n)$, namely 
$U_k = W_1^{v_k(1)}W_2^{v_k(2)}\dots$.  The 
unitary operators $U_1,U_2,\dots$ satisfy 
$$
U_iU_j=(-1)^{c_{ij}}U_jU_i \qquad i,j=1,2,\dots
\tag{2.9}
$$ 
and, after multiplication by 
suitable phase factors, 
$U_1, U_2,\dots$ becomes a spin system which generates 
the Clifford algebra $\Cal C$.  
\endexample

\subhead{3.  The Universal $C^*$-algebra}
\endsubhead

The purpose of this section is to prove

\proclaim{Theorem B}
Let $p=2,3,\dots$ be a prime and let $u_1,u_2,\dots$ be a 
universal sequence 
of unitary operators satisfying $u_k^p=\bold 1$ for all $k$ 
and the commutation relations (1.2).  Let 
$\omega: \Gamma\times\Gamma\to\Bbb Z_p$ be the skew-symmetric 
form (1.3) and let $n=2r$ be its rank, $r=1,2,\dots,\infty$.  

Then $C^*(u_1,u_2,\dots)$ is isomorphic to $C(X)\otimes \Cal B$, 
where $X$ is a totally  
disconnected compact metrizable space, and where 
$\Cal B=M_{p^r}(\Bbb C)$ if $r$ is 
finite and is a UHF algebra of 
type $p^\infty$ if $r=\infty$.  

The center $C(X)\otimes\bold 1$ is the closed linear 
span of the set of words $\{w_x: x\in \ker\omega\}$.  
$C^*(u_1,u_2,\dots)$ is 
simple iff its center is trivial iff $\omega$ is a symplectic form.  
\endproclaim

\remark{Remarks}
Since every quotient of $C(X)$ for $X$ a compact 
totally disconnected metrizable space is of 
the form $C(Y)$ for $Y$ of the same type, it follows that 
any sequence of unitary operators $U_1,U_2,\dots$ that satisfies 
$U_k^p=\bold 1$ and the relations (1.2), whether it is universal 
or not, must generate a $C^*$-algebra of the same general type 
$C(Y)\otimes \Cal B$ as the universal one $C(X)\otimes \Cal B$.  
If $\{U_1,U_2,\dots\}$ is irreducible and $\omega$ is 
of infinite rank, then $X$ reduces to a point and 
$C^*(U_1, U_2,\dots)$ 
is a UHF algebra of type $p^\infty$.  
\endremark

Before giving the proof of Theorem B, we require
two elementary results.  

\proclaim{Lemma 3.1}
Let $\Cal A$ be a unital $C^*$-algebra which is 
generated by two mutually  commuting 
unital $C^*$-subalgebras $\Cal Z$, $\Cal B$ with the 
properties
\roster
\item"{(i)}" $\Cal Z\cong C(X)$ is commutative, and 
\item"{(ii)}" $\Cal B$ is a UHF algebra.  
\endroster
Then $\Cal Z$ is the center of $\Cal A$ and 
$\Cal A\cong C(X)\otimes\Cal B$.  
\endproclaim

\demo{proof of Lemma 3.1}
The proof is straightforward and we merely sketch 
the argument.  Suppose first that the subalgebra $\Cal B$ 
is finite dimensional, hence isomorphic to the matrix 
algebra $M_n(\Bbb C)$ for some $n=1,2,\dots$.  Pick 
a set of matrix units $e_{ij}$, $1\leq i,j\leq n$ for 
$\Cal B$.  Thus 
$e_{ij}e_{kl}=\delta_{jk}e_{il}$, $e_{ij}^*=e_{ji}$, 
and $e_{11}+\dots+e_{nn}=\bold 1$.  Using these 
relations and the fact that the elements of $\Cal
Z$ commute with the $e_{ij}$ one finds that for 
arbitrary $z_{ij}\in\Cal Z$, $1\leq i,j\leq n$,
$$
\sum_{i,j=1}^n z_{ij}e_{ij}=0 \implies z_{ij}=0, 
\quad {\text{for all }} 1\leq i,j\leq n.  
$$
Thus if we consider $\Cal Z\otimes \Cal B$ to be 
the $C^*$-algebra $M_n(\Cal Z)$ then the preceding 
observation shows that the natural $*$-homomorphism 
$\pi: M_n(\Cal Z)\to \Cal A$ defined by 
$$
\pi((z_{ij}))=\sum_{i.j=1}^n z_{ij}e_{ij}
$$
is injective; it also has dense range, hence it is a 
$*$-isomorphism which carries the center of 
$M_n(\Cal Z)$ onto $\Cal Z$.  

In the general case, $\Cal B$ is the norm closure 
of an increasing sequence of algebras $\Cal B_n$ of 
the above type.  The preceding argument shows that 
the natural surjective $*$-homomorphism 
$\pi: \Cal Z\otimes\Cal B\to\Cal A$ restricts 
to an isometric $*$-homomorphism on 
each $\Cal Z\otimes\Cal B_n$, 
hence it is an isometric 
$*$-isomorphism.  
\qedd\enddemo

\proclaim{Lemma 3.2}
Let $p$ be a positive integer and let $V$ and $W$ be unitary 
operators in some $C^*$-algebra satisfying $V^p=W^p=\bold 1$
and $VW=\zeta WV$ 
where $\zeta=e^{2\pi i/p}$.  Then $C^*(V,W)\cong M_p(\Bbb C)$.  
\endproclaim

\demo{proof of Lemma 3.2}
Since $W^p=\bold 1$, the spectrum $\sigma(W)$ of $W$ is contained 
in the set of $p$th roots of unity, and because 
$VWV^{-1}=\zeta W$, $\sigma(W)$ is invariant under multiplication
by $\zeta$.  Hence $\sigma(W)=\{1,\zeta,\zeta^2,\dots,\zeta^{p-1}\}$.  
Letting $P_k$ be the spectral projection of $W$ corresponding to 
the eigenvalue $\zeta^k$, $k=0,1,\dots,p-1$, the commutation 
relation $V^iWV^{-i}=\zeta^{i}W^j$ implies that 
$V^iP_j=P_{i+j}V^i$, where the sum $i+j$ is interpreted modulo 
$p$.  Together with $V^p=\bold 1$, this implies that 
the operators $e_{ij}=V^{i-j}P_j$, $0\leq i,j\leq p-1$,
are a set of $p\times p$ matrix units which have
$C^*(V,W)$ as their linear span.  
\qedd\enddemo

\demo{proof of Theorem B}
Fix a universal sequence $u_1,u_2,\dots$ as above and 
let $\Cal Z$ be the closed linear span of 
the words of the form $w_x=u_1^{x_1}u_2^{x_2}\dots$ 
where $x\in\ker\omega$.  Notice that because 
$w_xw_y=(-1)^{\omega(x,y)}w_yw_x$, it follows that 
every word $w_x$ with $x\in\ker\omega$ belongs to 
the center of $C^*(u_1,u_2,\dots)$.  
Note that for each $x\in\Gamma$ 
we can choose a scalar $\lambda_x\in\Bbb T$ 
with the property that $(\lambda_xw_x)^p=\bold 1$.  
It is possible to do this because the relation 
$w_sw_t=\zeta^{Q(s,t)}w_{s+t}$ implies that 
$w_x^p$ is a scalar multiple of $w_{px}=w_0=\bold 1$.  
One can specify $\lambda_x$ explicitly, but it is 
not necessary to do so.  Thus $\Cal Z$ is a commutative 
AF algebra isomorphic to $C(X)$ for $X$ a compact 
metrizable totally disconnected space.  
Because of Lemma 3.1, it is enough to show that 
there is a UHF algebra $\Cal B\subseteq C^*(u_1,u_2,\dots)$ 
of the asserted type such that 
$C^*(u_1,u_2,\dots)$ is generated by 
$\Cal Z\cup\Cal B$.  

By Theorem A, $\Gamma$ decomposes into a direct sum of 
vector spaces $\Gamma=\ker\omega\oplus L$, where the restriction 
of $\omega$ to $L\times L$ is a symplectic form, and where 
$\dim L$ is the rank of $\omega$.  Since $L$ is a vector 
space, the relation $w_xw_y=\zeta^{Q(x,y)}w_{x+y}$ implies 
that $\Cal B=\overline{\text{span}}\{w_x: x\in L\}$
is a $C^*$-subalgebra of $\Cal A$.  Moreover, since 
$\Gamma=\ker\omega + L$, the set of products of 
words of the form 
$w_xw_y=\zeta^{Q(x,y)}w_{x+y}$, $x\in\ker\omega$, $y\in L$
have $\Cal A$ as their closed linear span.  Thus 
$\Cal Z\cup\Cal B$ generates $\Cal A$.  

It remains to show that $\Cal B$ is a UHF algebra of the 
asserted type. 
Suppose first that $\dim L=2r$ is
finite.  By  Theorem $A$ we can find a symplectic basis 
$e_1,\dots,e_r,f_1,\dots,f_r$ for the symplectic 
vector space $(L, \omega_L)$ obtained by restricting
$\omega$ to $L$.  
Consider  the operators $V_1,\dots,V_r,W_1,\dots,W_r$
defined by 
$$
V_k = \lambda_{e_k}w_{e_k}, \quad 
W_k = \lambda_{f_k}w_{f_k},  \tag{3.1}
$$
where the scalars $\lambda_x$ are as above.  
Every $x\in L$ is a linear combination of elements 
of $e_1,\dots,e_r,f_1,\dots,f_r$, hence the set of 
all products 
$V_1^{m_1}\dots V_r^{m_r}W_1^{n_1}\dots,W_r^{n_r}$, 
$m_1,\dots,m_r,n_1,\dots,n_r=0,1,\dots,p-1$, spans 
$\Cal B$.  
We have already arranged that $V_k^p=W_k^p=\bold 1$ for 
every $k$.  Note that for all $i,j=1,\dots,r$
$$
V_iV_j=V_jV_i, \quad
W_iW_j=W_jW_i, \quad
V_iW_j=\zeta^{\delta_{ij}}W_jV_i,   \tag{3.2}
$$  
$\delta_{ij}$ denoting the Kronecker delta.  
Indeed, these relations are immediate consequences 
of the 
basic formula $w_xw_y=\zeta^{\omega(x,y)}w_yw_x$ 
and the fact that $\{e_i,f_j\}$ is a
symplectic  set for $\omega$.  It follows from (3.2) that the
$C^*$-algebras 
$C^*(V_i,W_i)$ and $C^*(V_j,W_j)$ commute for 
$i\neq j$; and by Lemma 3.2 each 
$C^*(V_k,W_k)$ is isomorphic to $M_p(\Bbb C)$.  
Thus $\Cal B$ is isomorphic to a tensor  
product of $r$ compies of $M_p(\Bbb C)$.  

If $\dim L$ is infinite, then another application 
of Theorem A provides an infinite symplectic basis 
$e_1, e_2,\dots,f_1, f_2,\dots$ for $L$.  We define 
$V_1,V_2,\dots, W_1, W_2,\dots$ by (3.1) as before, 
and these operators satisfy (3.2).  In this case, 
the $C^*$-algebra 
$\Cal B$ generated by $V_i,W_j$ commutes with $\Cal Z$, 
and is generated by an increasing sequence of 
subalgebras $\Cal B_1\subseteq\Cal B_2\subseteq\dots$ 
$$
\Cal B_n=C^*(V_1,\dots,V_n,W_1,\dots,W_n), \qquad n=1,2,\dots.
$$ 
The preceding paragraph shows that $\Cal B_n$ is 
isomorphic 
$M_{p^n}(\Bbb C)$.  Hence $\Cal B$ is a UHF algebra 
of type $p^\infty$.  

The assertions of the third paragraph of 
Theorem B are obvious 
consequences of what has already been proved.  
\qedd\enddemo

\subhead{4.  Irreducible Spin Systems}
\endsubhead

Let $C=(c_{ij})$ be a matrix of zeros and ones, fixed throughout
the remainder of this section; in order to rule out 
the degeneracies described in the introduction, we also assume
that $(c_{ij})$ is of infinite rank.  Thus, 
the $\Bbb Z_2$-valued bilinear form 
$$
\omega(x,y)=\sum_{p,q=1}^\infty c_{pq}x_qy_p, \qquad x,y\in \Gamma
\tag{4.1}
$$
associated with $C=(c_{ij})$ has the property that 
$\Gamma/\ker\omega$ is infinite dimensional, 
$\ker\omega$ being the linear subspace 
$\{x\in\Gamma: \omega(x,\Gamma)=\{0\}\}\subseteq\Gamma$.

The purpose of this section is 
to classify the irreducible spin systems associated with 
$C$. 
Thus we consider irreducible spin
systems $\bar U=(U_1,U_2,\dots)$ acting on 
an infinite dimensional 
Hilbert space $H$, satisfying 
$$
U_iU_j=(-1)^{c_{ij}}U_jU_i,\qquad i,j=1,2\dots .\tag{4.2}
$$
Theorem B implies that $C^*(U_1,U_2,\dots)$ is the CAR algebra, 
and since the CAR algebra is a simple $C^*$-algebra not of 
type I, there can be no meaningful classification of 
such sequences $\bar U$ up to unitary equivalence.  
The equivalence relation that is 
appropriate for irreducible spin systems is 
Voiculescu's notion of approximate unitary equivalence \cite{A2}.  
Two spin systems $\bar U$ and $\bar V$, acting 
on infinite dimensional Hilbert spaces $H$ and $K$, respectively, 
are said to be equivalent (written $\bar U\sim\bar V$) 
if there is a sequence 
of unitary operators $W_1, W_2,\dots: H\to K$ such that 
$$
\lim_{n\to\infty}\|W_nU_kW_n^{-1}-V_k\|=0\qquad k=1,2,\dots.  
$$

We first introduce an invariant for 
irreducible spin systems $\bar U$.  
For every $x\in\Gamma$ there is a word 
$$
W_x=U_1^{x_1}U_2^{x_2}\dots,\qquad x\in\Gamma, 
$$
and we have $W_xW_y=(-1)^{Q(x,y)}W_{x+y}$ for 
all $x,y\in\Gamma$, where 
$Q:\Gamma\times \Gamma\to\Bbb Z_2$ is the 
bilinear form (1.5).
If $x\in\ker\omega$ then by (1.4) $W_x$ commutes with 
all words, and by irreducibility it must be 
a scalar multiple of the identity
$$
W_x=f(x)\bold 1, \qquad x\in\ker\omega.  
$$
This defines a function $f: \ker\omega\to\Bbb T$ 
satisfying the functional equation 
$$
f(x)f(y)=(-1)^{Q(x,y)}f(x+y), \qquad x,y\in\ker\omega.  
\tag{4.3}
$$
$f$ is called the {\bf standard invariant} associated 
with the irreducible spin system $\bar U$.  
Notice that (4.3) implies that $f(0)=1$.  
Since $f(x)^2=(-1)^{Q(x,x)}f(2x)=(-1)^{Q(x,x)}f(0)=\pm 1$, 
it follows that $f$ must take values in the 
multiplicative group of fourth roots of unity, 
$$
f(x)^4 = 1, \qquad x\in\ker\omega.  
$$

\proclaim{Proposition 4.1}
Let $\bar U=(U_1,U_2,\dots)$ and 
$\bar U^\prime=(U_1^\prime,U_2^\prime,\dots)$ be two irreducible 
spin systems on Hilbert spaces $H$, $H^\prime$ 
which satisfy the relations (4.2), and 
let $\pi$, $\pi^\prime$ be the representations 
of the universal $C^*$-algebra $\Cal A_C=C^*(u_1,u_2,\dots)$ 
defined 
by $\pi(u_k)=U_k$, $\pi^\prime(u_k)=U_k^\prime$, 
$k=1,2,\dots$.  The following are equivalent.  
\roster
\item"{(i)}"
$\bar U\sim\bar U^\prime$.  
\item"{(ii)}"
$\ker \pi=\ker\pi^\prime$.  
\item"{(iii)}"
$\bar U$ and $\bar U^\prime$ have the same standard invariant.  
\item"{(iv)}"
For every $n=1,2,\dots$, there is a unitary operator 
$W_n: H\to H^\prime$ such that 
$$
W_n U_k W_n^{-1}=U_k^\prime, \qquad k=1,2,\dots,n.  
$$
\endroster
\endproclaim

\demo{proof}
The implications (i)$\implies$ (ii)$\implies$ (iii) and 
(iv)$\implies$ (i) are straightforward.  
We prove (iii)$\implies$ (iv).   
For that we will make use of the following elementary result.
The proof, 
a straightforward exercise in elementary 
multiplicity theory (see \cite{A1}), is omitted.  

\proclaim{Lemma 4.1}
Let $\Cal B$ be a finite dimensional $C^*$-algebra and 
let $\pi_1$, $\pi_2$ be two 
faithful nondegenerate representations 
of $\Cal B$ on Hilbert spaces $H_1$, $H_2$ such that 
$\pi_j(\Cal B)\cap\Cal K_j=\{0\}$ for $j=1,2$, 
$\Cal K_j$ denoting 
the compact operators on $H_j$.  Then 
$\pi_1$ and $\pi_2$ are unitarily equivalent.  
\endproclaim

Let $f$ and $f^\prime$ be the respective standard invariants 
of $\bar U$ and $\bar U^\prime$.  
Assuming that $f=f^\prime$ as in (iii),  
we have to verify (iv), and by replacing the spin 
system $\bar U^\prime$ with a unitarily equivalent 
one, we may assume that both $\bar U$ and 
$\bar U^\prime$ act on the same Hilbert space.  
Consider the words $w_x=u_1^{x_1}u_2^{x_2}\dots$ 
in $\Cal A_C$ corresponding to elements $x\in\ker\omega$.  
Since $f=f^\prime$ we have 
$$
\pi(w_x)=U_1^{x_1}U_2^{x_2}\dots=f(x)\bold 1=
f^\prime(x)\bold 1=\pi^\prime(w_x), \qquad x\in\ker\omega.  
$$
Since by Theorem B 
the central words of this type have the center 
$\Cal Z$ of $\Cal A_C$ as their closed linear span, 
it follows that 
$\pi\restriction_\Cal Z=\pi^\prime\restriction_\Cal Z$.  
Theorem B also implies that $\Cal A_C$ is isomorphic 
to $C(X)\otimes\Cal C$ where $\Cal C$ is the CAR 
algebra,  hence any two irreducible representations that 
agree on the center must have the same kernel (corresponding 
to some point $p\in X$).  Hence $\ker \pi=\ker \pi^\prime$.  

It follows that for every operator $A\in\Cal A_C$ 
we have $\|\pi(A)\|=\|\pi^\prime(A)\|$, hence
there is a unique $*$-isomorphism 
$\alpha: C^*(U_1,U_2,\dots)\to C^*(U_1^\prime,U_2^\prime,\dots)$
such that $\alpha\circ\pi=\pi^\prime$.  
Both of these are simple unital $C^*$-algebras  (they are 
isomorphic to the CAR algebra), and hence contain no 
nonzero compact operators.  Noting that the restriction
of $\alpha$ to $C^*(U_1,\dots,U_n)$ is a $*$-isomorphism 
onto $C^*(U_1^\prime,\dots,U_n^\prime)$ which carries 
the $n$-tuple of operators $(U_1,\dots,U_n)$ to 
$(U_1^\prime,\dots,U_n^\prime)$, an application of 
Lemma 4.1 implies that each of these restrictions 
is implemented by a unitary operator $W_n\in\Cal B(H)$, 
and (iv) follows.  
\qedd
\enddemo

In conclusion, we describe how the irreducible spin 
systems associated with a commutation matrix $C$ 
can be described and classified in terms 
of any one of them.  For any irreducible spin
system $\bar U=(U_1,U_2,\dots)$, we consider the 
spin systems that can be obtained from it by 
changing phases as follows.  For every sequence 
of numbers $\gamma=(\gamma_1,\gamma_2,\dots)$ 
in $\{0,1\}=\Bbb Z_2$ consider the 
sequence of unitary operators
$$
\bar U^\gamma=((-1)^{\gamma_1}U_1,(-1)^{\gamma_2}U_2,\dots).  
$$
It is clear that 
$\bar U^\gamma$ is an irreducible spin system 
satisfying the same commutation relations as $\bar U$.  
We now show that these ``phase shifted" versions 
of $\bar U$ provide all possible standard invariants.  

\proclaim{Lemma 4.2}
Let $\bar U$ be an irreducible spin system, let 
 $f: \ker\omega\to\Bbb T$ be its standard invariant, and 
let $g: \ker\omega\to\Bbb T$ be any function 
satifying the same functional equation (4.3)
$$
g(x)g(y)=(-1)^{Q(x,y)}g(x+y), \qquad x,y\in\ker\omega.  
$$  
Then there is a 
$\gamma=(\gamma_1,\gamma_2,\dots)\in\Bbb Z_2^\infty$ 
such that $g$ is the standard invariant of 
$\bar U^\gamma$.  
\endproclaim

\demo{proof}
For $\gamma\in\Bbb Z_2^\infty$, we can express 
the standard invariant $f^\gamma$ of $\bar U^\gamma$ in terms 
of the standard invariant $f$ of $\bar U$ as follows.  
For every $x\in\Gamma$ the word for 
$\bar U^\gamma$ is 
$$
(-1)^{\sum_k\gamma_kx_k}U_1^{x_1}U_2^{x_2}\dots, 
$$
hence for $x\in\ker\omega$ we have 
$$
f^\gamma(x)=(-1)^{\sum_k\gamma_kx_k}f(x).  \tag{4.4}
$$

Now both $g$ and $f$ satisfy (4.3), hence 
the function $h: \ker\omega\to\Bbb T$ defined by 
$h(x)=g(x)/f(x)$ satisfies 
$$
h(x+y)=h(x)h(y), \qquad x,y\in\ker\omega.  
$$
Notice too that since $x+x=0$ for all $x\in\ker\omega$ 
we have 
$h(x)^2=h(x)h(x)=h(x+x)=h(0)=1$.  It follows that 
$h(x)=\pm 1$ for all $x\in\ker\omega$.  Thus there 
is a unique function $\theta: \ker\omega\to\{0,1\}=\Bbb Z_2$ satisfying 
$$
g(x)/f(x)=h(x)=(-1)^{\theta(x)}, \qquad x\in\ker\omega, \tag{4.5}
$$
and we have $\theta(x+y)=\theta(x)+\theta(y)$ relative to 
the addition in the field $\Bbb Z_2$ because 
$h(x+y)=h(x)h(y)$ for $x,y\in\ker\omega$.  

We may consider $\theta: \ker\omega\to\Bbb Z_2$ as a linear  
functional defined on the vector space 
$\ker\omega\subseteq\Gamma$.  A familiar argument implies 
that a linear functional defined on a subspace of a vector 
space can be extended to a linear functional defined on 
the entire space.  Thus we may find a function 
$\tilde\theta: \Gamma\to\Bbb Z_2$ such that 
$\tilde\theta(x+y)=\tilde\theta(x)+\tilde\theta(y)$ 
for all $x,y\in\Gamma$ and which restricts to $\theta$ 
on $\ker\omega$.  Letting $u_1, u_2,\dots$ be the usual 
basis of unit vectors for $\Gamma$, $u_k(j)=\delta_{kj}$, 
we define $\gamma=(\gamma_1,\gamma_2,\dots)\in\Bbb Z_2^\infty$ 
by $\gamma_k=\tilde\theta(u_k)$, $k=1,2,\dots$.  
For every $x=(x_1,x_2,\dots)\in\Gamma$ we have 
$\tilde\theta(x)=\sum_{k=1}^\infty\tilde\theta(u_k)x_k=
\sum_{k=1}^\infty \gamma_kx_k$.  
Substituting the latter into (4.5) we find that 
$$
g(x)=(-1)^{\theta(x)}f(x)=
(-1)^{\sum_{k=1}^\infty \gamma_k x_k}f(x),\qquad
x\in\ker\omega.  
$$
By (4.4), this is the standard invariant $f^\gamma$ of 
$\bar U^\gamma$.  
\qedd\enddemo

\proclaim{Theorem C}Let $\bar U$ be any irreducible 
spin system satisfying the commutation relations (1.1) and 
let $\omega$ be the skew-symmetric form (1.3).  
Every irreducible spin sytem satisying the same 
commutation relations is equivalent 
to $\bar U^\gamma$ for some $\gamma\in\Bbb Z_2^\infty$.  
Given two sequences $\gamma$, $\gamma^\prime$ in 
$\Bbb Z_2^\infty$, the spin systems 
$\bar U^{\gamma}$ and $\bar U^{\gamma^\prime}$ 
are equivalent 
if and only if $\gamma$ and $\gamma^\prime$ define the same linear 
functional on $\ker\omega$,  
$$
\sum_{k=1}^\infty \gamma_kx_k=\sum_{k=1}^\infty \gamma^\prime_kx_k,
\qquad x\in\ker\omega.  
$$

In particular, if $\ker\omega$ is of
finite dimension $d$ as a vector space over $\Bbb Z_2$, 
then there are exactly $2^d$ 
equivalence classes of irreducible spin systems 
associated with $C$.  If $\ker\omega$ is infinite 
dimensional then the set of distinct equivalence 
classes of irreducible spin systems has the cardinality 
of the continuum $2^{\aleph_0}$.  
\endproclaim

\demo{proof}
Fix an irreducible spin system $\bar U$ as above, and 
let $\gamma$ and $\gamma^\prime$ be two sequences in 
$\Bbb Z^\infty$.  We show first that 
$\bar U^\gamma\sim\bar U^{\gamma^\prime}$ $\iff$ 
$$
\sum_{k=1}^\infty (\gamma_k)x_k=
\sum_{k=1}^\infty\gamma_k^\prime x_k, 
\qquad x\in\ker\omega.  \tag{4.6}
$$
Indeed, lettting 
$f^\gamma$ and $f^{\gamma^\prime}$ be the standard 
invariants for $\bar U^\gamma$ and $\bar U^{\gamma^\prime}$, 
we see from (4.4) that 
$$
f^\gamma(x)=(-1)^{\sum_k\gamma_kx_k}f(x), \quad 
f^{\gamma^\prime}(x) = 
(-1)^{\sum_k\gamma^\prime_kx_k}f(x), \qquad x\in\ker\omega
$$
and hence $f^\gamma=f^{\gamma^\prime}$ iff (4.6) holds. 
 By the characterization (iii) of 
Proposition (4.1), that is equivalent to 
$\bar U^\gamma\sim\bar U^{\gamma^\prime}$.

Now let $V=(V_1,V_2,\dots)$ be an arbitrary 
irreducible spin system associated with $C$, and 
let $g: \ker\omega\to\Bbb T$ be its standard invariant.  
Lemma 4.2 implies that there is a 
$\gamma\in\Bbb Z_2^\infty$ such that $g=f^\gamma$, 
and by part (iii) of Proposition 4.1 we conclude 
that $\bar V\sim\bar U^\gamma$.

It remains only to establish the results on cardinality, 
and in view of what has been proved, we simply 
have to count the distinct functions $g: \ker\omega\to\Bbb T$ 
that satisfy the functional equation (4.3).  Letting $f$ be 
the standard invariant of $\bar U$, the proof of 
Lemma 2 shows that every such $g$ is 
obtained from it by way of 
$$
g(x) = (-1)^{\theta(x)}f(x), \qquad x\in\ker\omega, 
$$
where $\theta: \ker\omega\to\Bbb Z_2$ is a 
(necessarily unique) linear functional.  
Thus the set of all such $g$ is in bijective correspondence 
with the set of linear functionals on $\ker\omega$.  
If $\ker\omega$ is of finite dimension $d$ then by choosing a basis 
$e_1,\dots,e_d$ for $\ker\omega$ we find that the set of all such 
$\theta$ is in bijective correspondence with the set of 
all functions from $\{e_1,\dots,e_d\}$ to $\Bbb Z_2$, and 
the cardinality of that set is $2^d$.  

If $\ker\omega$ is infinite dimensional, then since it is 
a countably generated vector space 
it has a basis $\{e_1, e_2,\dots\}$.  As in the 
preceding paragraph, the set of all standard invariants is in 
bijective correspondence with the set of all linear functionals 
on $\ker\omega$, which in turn corresponds bijectively with 
the set of all functions from $\{e_1,e_2,\dots\}$ to 
$\Bbb Z_2$, a set of cardinality $2^{\aleph_0}$.  
\qedd\enddemo

When the commutation matrix is symplectic 
one has the following uniqueness:  

\proclaim{Corollary}
Let $C=(c_{ij})$ be an infinite matrix of zeros and ones which 
is skew-symmetric and nondegenerate.  Then any two irreducible 
spin systems satisfying the commutation relations 
$U_iU_j=(-1)^{c_{ij}}U_jU_i$ are approximately unitarily 
equivalent.  
\endproclaim

\remark{Remark}
In such cases the $C^*$-algebra $\Cal A_C$ associated 
with $C$ is the CAR algebra, and is therefore simple 
not of type $I$.  In view of Proposition 4.1, the Corollary 
remains valid {\it verbatim} if one deletes the irreducibility 
hypothesis.  
\endremark


\enddocument